\documentclass[12pt]{amsart}
\input{xypic}
\usepackage{amsmath}
\usepackage{amssymb}
\usepackage{amscd}
\usepackage{latexsym}
\usepackage{enumerate}

\newtheorem{theorem}{Theorem}[section]
\newtheorem{proposition}[theorem]{Proposition}

\newtheorem{lemma}[theorem]{Lemma}

\newtheorem{remark}[theorem]{Remark}

\newtheorem{definition}[theorem]{Definition}

\newcommand{\M}{{\mathcal M}}

\newcommand{\Z}{{\mathbb Z}}
\newcommand{\N}{{\mathbb N}}

\newcommand{\dst}{\displaystyle}

\newcommand{\im}{{\rm Im}}

\newcommand{\downarrowright}[1]{\downarrow
\rlap{\raise0.1cm\hbox{$\scriptstyle{#1}$}}}

\newcommand{\downarrowleft}[1]{\rlap{\kern-0.2cm
\raise0.1cm\hbox{$\scriptstyle{#1}$}}\downarrow}

\newcommand{\uparrowright}[1]{\uparrow
\rlap{\lower0.1cm\hbox{$\scriptstyle{#1}$}}}

\newcommand{\uparrowleft}[1]{\rlap{\kern-0.2cm
\lower0.1cm\hbox{$\scriptstyle{#1}$}}\uparrow}

\newcommand{\la}{\leftarrow}

\newcommand{\ra}{\rightarrow}
\newcommand{\lra}{\longrightarrow}

\newcommand{\epi}{\mbox{$\to$\hspace{-0.35cm}$\to$}}

\newcommand{\mono}{\hookrightarrow}

\def\rmono{\rto|<\hole|<<\ahook}

\def\umono{\ar@{_{(}->}[u]}
\def\uumono{\ar@{_{(}->}[uu]}

\def\lmono{\ar@{_{(}->}[l]}
\def\llmono{\ar@{_{(}->}[ll]}

\begin{document}
\title{Homology fibrations and ``group-completion'' revisited}

\author{Wolfgang Pitsch and J\'er\^ome Scherer}

\date{}



\begin{abstract}
We give a proof of the Jardine-Tillmann generalized group
completion theorem. It is much in the spirit of the original
homology fibration approach by McDuff and Segal, but follows a
modern treatment of homotopy colimits, using as little simplicial
technology as possible. We compare simplicial and topological
definitions of homology fibrations.
\end{abstract}

\maketitle


\section*{Introduction}
The group completion of a topological monoid $M$ is the loop space
$\Omega BM$ and a group completion theorem is originally a
statement about the relation between the homology of $M$ and that
of $\Omega BM$. In the appendix of \cite{MR95a:14023} D.~Quillen
considers a simplicial monoid $M$. His main theorem is that under
certain conditions the homology of the group completion of $M$ can
be computed by inverting $\pi_0 M$ in the homology of~$M$. A
similar result can be found in May's
\cite[Theorem~15.1]{MR51:6806}. In this paper we focus on a more
topological kind of group completion theorem, the question being
how to construct $\Omega BM$ out of $M$. Our starting point is
McDuff's and Segal's theorem, as it can be found
in~\cite[Proposition 2]{MR53:6547} (a good account on the subject
is Adams' book on infinite loop spaces
\cite[Chapter~3]{MR80d:55001}).

\medskip

\noindent{\bf Theorem}\,
{\it Let $M$ be a topological monoid acting on a space $X$ by
homology equivalences. Then the map $\pi: EM \times_M X \ra BM$
from the Borel construction to the classifying space of $M$ is a
homology fibration with fibre $X$.}

\medskip

The standard application is as follows. Let $M$ be a homotopy
commutative topological monoid with $\pi_0 M \cong \N$. Choose a
point $m$ in the component of $1$ and form the telescope $M_\infty
= Tel(M \stackrel{\cdot m}{\lra} M \stackrel{\cdot m}{\lra}
\dots)$. The action of $M$ by left multiplication on $M_\infty$ is
by homology equivalences because $M$ is homotopy commutative.
Hence we obtain:

\medskip

\noindent{\bf Corollary}\,
{\it Let $M$ be a homotopy commutative topological monoid. Then there is
a homology equivalence $M_\infty \ra \Omega BM$. Moreover, when $\pi_1 M_\infty$
is perfect, $\Omega BM \simeq M_\infty^+$.
}

\medskip

Taking for example $M$ to be the disjoint union $\coprod
B\Sigma_n$ of classifying spaces of the symmetric groups, the
Barrat-Priddy-Quillen Theorem states that $B \Sigma_\infty^+$ is
the infinite loop space $QS^0$, \cite{MR47:3489}. Likewise, taking
$M$ to be $\coprod BGL_n(R)$ one gets back Quillen's definition of
the algebraic $K$-theory of a ring $R$, \cite{MR43:5525}.

\medskip

Simplicial versions of the group completion theorem started
appearing at the end of the eighties. I.~Moerdijk provides a
homological statement in \cite[Corollary~3.1]{MR91j:18021} and
J.F.~Jardine the analogue of the above theorem in \cite[Theorem
4.2]{MR90j:55029}, which he calls the ``strong form of the Group
Completion Theorem". More recently U.~Tillmann introduced a
``multiple object case" in her celebrated work on the stable
mapping class group (\cite[Theorem 3.2]{MR99k:57036}). In this
context the Borel construction is replaced by a bisimplicial
version, i.e. the realization of a certain simplicial space. Let
$\M$ be a simplicial category and $F: \M^{op} \ra Spaces$ a
contravariant diagram. There is always a natural transformation to
the trivial diagram. Taking the bisimplicial Borel constructions
yields a map $\pi_\M: E_\M F \ra B\M$, analogous to the map $\pi$
in the classical theorem.

\medskip

{\bf Theorem~\ref{groupcompletion}.}
{\it Let $\M$ be a simplicial category and $F: \M^{op} \ra Spaces$
a contravariant diagram. Assume that any morphism $f: i \ra j$
induces an isomorphism in integral homology $H_*(F(j); \Z) \ra
H_*(F(i); \Z)$. Then, for each object $i \in \M$, the map $F(i)
\ra Fib_i(\pi_\M)$ to the homotopy fibre of $\pi_\M$ over $i$ is a
homology equivalence.}

\medskip

We offer in this paper a proof which uses as little simplicial
technology as possible. The main ingredient is a rather classical
result about comparing the fibre of the realization with the
realization of the fibres, an idea already used by McDuff and
Segal in their proof of the classical group completion theorem. Of
course we do not avoid simplicial spaces, the theorem after all is
about delooping a simplicial classifying space. We work however
more in the spirit of the modern theory homotopy colimits. One
very powerful tool in this setting is to decompose a space as a
diagram over its simplices. The advantage of this approach is that
one gets a more geometric feeling about the constructions
performed (such as the bisimplicial Borel construction). We also
use a simplicial notion of homology fibrations (preimages of
simplices have the same integral homology as the homotopy fibre).
In the last section we compare this concept to that of classical
homology fibration in the category of topological spaces and prove
they coincide.

In this paper space means simplicial set and we write $Spaces$ for
the category of spaces.

\medskip

\noindent {\bf Acknowledgements}: We would like to thank
J.B.~Bost, F.~Loeser, and I.~Madsen for organizing the week about
the mapping class group in Luminy (January 2002), and
W.~Chach\'olski for helpful comments. The second author was
introduced to simplicial technology by W.~Chach\'olski, so
probably he and W.~Dwyer have a similar proof of the group
completion theorem in a drawer.

\section{Homology fibrations}
\label{homologyfibrations}

Let $p: E \ra B$ be a map of spaces and $\sigma$ be an $n$-simplex in $B$.
We denote by $dp(\sigma)$ the pull-back of the diagram
$\Delta[n] \stackrel{\sigma}{\lra} B \stackrel{p}{\longleftarrow} E$. This is the ``preimage"
of the simplex in $E$ and yields a functor $dp: \Delta B \ra Spaces$ from the simplex category
of the base space (this category is defined for example in \cite[p.182]{dror:book}, see also
\cite[Definition 6.1]{ChSc}). It allows to decompose the map $p$ as a diagram over $\Delta B$,
as one has $E \simeq hocolim_{\Delta B} dp$ and $B \simeq hocolim_{\Delta B} \Delta[n]$.

We will also need a slight generalization of $dp$, replacing a simplex by any space $K$.
For a map $f: K \ra B$, define $dp(f)$ to be the pull-back of $f$ along $p$.

\begin{definition}
{\rm A map of spaces $p: E \ra B$ is a {\em homology fibration} if the
natural map $dp(\sigma) \ra Fib_\sigma(p)$ to the homotopy fibre of $p$
over the component of $\sigma$ is a homology equivalence for any simplex
$\sigma \in B$. It is a {\em weak homology fibration} if for any simplex
$\sigma \in B$ and any simplicial operation $\theta$ we have a homology equivalence
$dp(\sigma) \ra dp(\theta \sigma)$.}
\end{definition}

The aim of this section is to prove that a weak homology fibration is
actually a homology fibration. This part of the paper replaces Segal and
McDuff's work on locally contractible paracompact spaces.

\begin{lemma}{\rm \cite[Proposition 6]{MR53:6547}}
\label{contractiblebase}
Let $p: E \ra B$ be a weak homology fibration with $B$ contractible. Then
$p$ is a homology fibration.
\end{lemma}

\begin{proof}
The category $\Delta B$ is contractible since $B \simeq hocolim_{\Delta B} *
= N(\Delta B)$. So $E$ is equivalent to the homotopy colimit over a contractible
category of a diagram in which all maps are homology equivalences. This homotopy
colimit has the same homology type as any of the values $dp(\sigma)$ since it can
be computed (\cite{Amit}) by using only push-outs and telescopes of diagrams
consisting of homology equivalences. We conclude by the Mayer--Vietoris Theorem
and the fact that homology commutes with telescopes.
\end{proof}

\begin{proposition}
\label{pullbackhomologyfibration}
Let $p: E \ra B$ be a weak homology fibration and $f: B' \epi B$ a fibration.
The pull-back of $p$ along $f$ is another weak homology fibration $p': E' \ra B'$.
\end{proposition}

\begin{proof}
Let $\sigma'$ be a simplex in $B'$, $\sigma = f\sigma'$ its image in $B$ and
$\theta$ any simplicial operation. Then $dp(\sigma)$ has the same homology
type as $dp(\theta \sigma)$ by assumption. But
$dp'(\sigma') \simeq dp(\sigma)$ and $dp'(\theta\sigma') \simeq dp(\theta\sigma)$
since $p'$ was obtained as a pull-back.
\end{proof}

\begin{theorem}{\rm \cite[Proposition 5]{MR53:6547}}
\label{weakandstrong}
A weak homology fibration is a homology fibration.
\end{theorem}

\begin{proof}
Let $p: E \ra B$ be a weak homology fibration and choose $f: PB \epi B$ the
path space fibration. The above proposition applies,
so $p': Fib_\sigma (p) \ra PB$ is a weak homology fibration as well for
any simplex $\sigma$ in $B$.
Since $f$ is surjective, there exists a simplex $\sigma' \in PB$ such that
$f(\sigma') = \sigma$. Therefore $dp(\sigma) \simeq dp'(\sigma')$, which
has the same homology type as the homotopy fibre $Fib_\sigma (p)$ by
Lemma~\ref{contractiblebase}.
\end{proof}

\section{Realizations and fibres}
\label{realization}

Theorem~\ref{weakandstrong} will be used throughout this section. For
checking that a map is a homology fibration it suffices to check it is
a weak homology fibration.

\begin{lemma}
\label{preimage}
Consider a commutative square
\[\xymatrix{
E_{0}\rto\dto_{p_{0}} & E_{1}\dto^{p_{1}}\\
B_{0}\rto & B_{1}
}\]
where the vertical arrows are compatible homology fibrations in the sense that the map
$Fib_{v}(p_0) \ra Fib_{v} (p_1)$ is an integral homology equivalence for any vertex
$v \in B_0$. Then $dp_0(f) \ra dp_1(f)$ is an integral homology equivalence for any
map $f: K \ra B_0$. Moreover if both horizontal maps are cofibrations, then so is
$dp_0(f) \ra dp_1(f)$.
\end{lemma}

\begin{proof}
Notice first that if $\sigma$ is a simplex in $B_0$, then
$dp_0(\sigma) \ra dp_1(\sigma)$ is an integral
homology equivalence by our assumption on the homotopy fibres over vertices.
Likewise the preimages in $E_0$ and $E_1$ of a disjoint union of simplices
have the same integral homology type. We assume therefore that $K$ is connected.
Assume $K = L \cup_{\partial \Delta[n]} \Delta[n]$. By induction on the dimension
suppose that both $dp_0(f|_L) \ra dp_1(f|_L)$ and $dp_0(f|_{\partial \Delta[n]})
\ra dp_1(f|_{\partial \Delta[n]})$ are homology equivalences. We see that
the preimage of $\partial \Delta[n]$ is contained in that of $\Delta[n]$
so that
$$
dp_0(f) = colim \big( dp_0(f|_L) \la dp_0(f|_{\partial \Delta[n]})
\mono dp_0(f|_{\Delta[n]}) \big)
$$
is actually a homotopy push-out. Thus $dp_0(f) \ra dp_1(f)$ is a homotopy push-out
of homology equivalences.
\end{proof}

\medskip

We prove now that a push-out of homology fibrations is still a homology fibration. As everybody
knows a map can always be replaced by a fibration, so we must pay close attention to the
constructions we perform. We always use strict colimits, but for
diagrams where the colimit is weakly equivalent to the homotopy colimit.

\begin{proposition}
\label{pushouthomologyfibration}
Consider a natural transformation between push-out diagrams:
\[\xymatrix{
E\dto_{p}\ar @{}[r]|(0.46)= &
*{colim\hspace{2mm}\big(\hspace{-20pt}} & E_{1}\dto_{p_{1}}
&E_{0}\lto\rmono\dto_{p_{0}} & E_{2}\dto^{p_{2}} & *{\hspace{-20pt}\big)}\\
B \ar @{}[r]|(0.46)=& *{colim\hspace{2mm}\big(\hspace{-20pt}} &
B_{1} & B_{0}\lto\rmono & B_{2} & *{\hspace{-20pt}\big)}
}\]
such that $p_n: E_n \ra B_n$ is a homology fibration for $0 \leq n \leq 2$
and the right hand-side horizontal maps are cofibrations.
Assume that the map $Fib_{v}(p_0) \ra Fib_{v} (p_n)$
is an integral homology equivalence for any vertex $v \in B_0$ if $n = 1, 2$.
Then $p$ is a homology fibration as well. Moreover, if for some $0 \leq n \leq 2$,
$w$ is a vertex in $B_n$,
then $B_n \mono B$ induces a homology equivalence $Fib_{w} (p_n) \ra Fib_w(p)$.
\end{proposition}

\begin{proof}
Any simplex $\sigma$ in $B$ lies either in $B_1$ or
in $B_2$. Say it lies in $B_1$ (the other case is similar)
and consider the pull-back $K$ of $\Delta[n] \ra B_1 \la B_0$. Apply Lemma~\ref{preimage}
to the map $f: K \ra B_0$ to conclude
that $dp_0(f) \ra dp_2(f)$ is a homology equivalence, which is even
a cofibration. Hence the preimage $dp(\sigma)$ is the (homotopy) push-out
$\dst colim \big( dp_1(\sigma) \la dp_0(f) \mono dp_2(f) \big)$. The
homotopy push-out of a homology equivalence is again a
homology equivalence so that $dp(\sigma)$ has the same homology type
as $dp_1(\sigma)$. We conclude that $p$ is a weak homology fibration.
\end{proof}

\begin{proposition}
\label{telescopehomologyfibration}
Consider a natural transformation between telescope diagrams:
\[
\xymatrix{
E\dto_{f}\ar @{}[r]|(0.46)= &
*{colim\hspace{2mm}\big(\hspace{-20pt}} & E_{0}\dto_{p_{0}}
\rmono & E_{1}\rmono\dto_{p_{1}} & E_{2}\dto_{p_{2}} \rmono & \cdots &
*{\hspace{-20pt}\big)}\\
B \ar @{}[r]|(0.46)=& *{colim\hspace{2mm}\big(\hspace{-20pt}} &
B_{0}\rmono & B_{1}\rmono & B_{2}\rmono & \cdots &
*{\hspace{-20pt}\big)}
}\]
such that
$p_n: E_n \ra B_n$ is a homology fibration for any $n \geq 0$ and all
horizontal maps are cofibrations.
Assume that the map $Fib_{v}(p_n) \ra Fib_{v}(p_{n+1})$ is an
integral homology equivalence for any $n \geq 0$ and any vertex $v \in B_n$.
Then $p$ is a homology fibration as well.  Moreover, if $w$ is a vertex in $B_n$
for some $n \geq 0$,
then the inclusion $B_n \mono B$ induces a homology equivalence $Fib_w(p) \ra Fib_{w} (p_n)$.
\end{proposition}

\begin{proof}
As $B = \bigcup B_n$, any simplex $\sigma$ of $B$ lies in some $B_N$. The conclusion
follows since $dp(\sigma) = \bigcup_{n \geq N} dp_n(\sigma)$ has the same homology type
as $dp_N(\sigma)$.
\end{proof}

\medskip

Let $X_\bullet$ be a simplicial space. Recall that Segal's thick realization
$|| X_\bullet ||$ (\cite[Appendix A]{MR50:5782}) is defined by an inductive process.
We have $|| X_\bullet || = \bigcup_n || X_\bullet ||_n$ where $|| X_\bullet ||_0 = X_0$
and $|| X_\bullet ||_n$ is constructed from $|| X_\bullet ||_{n-1}$ by the following
push-out
$$
colim \big( || X_\bullet ||_{n-1} \la \partial\Delta[n] \times X_n \mono
\Delta[n] \times X_n \big)
$$
and the map $\partial\Delta[n] \times X_n \ra || X_\bullet ||_{n-1}$ is defined
using only the face maps. This thick realization can be seen as the homotopy colimit
of the diagram $X_\bullet$ over the subcategory of $\Delta^{op}$ generated by the
face morphisms.

\begin{theorem}{\rm \cite[Proposition 4]{MR53:6547}}
\label{simplicialhomologyfibration}
Let $p_\bullet: E_\bullet \ra B_\bullet$ be a map of simplicial spaces such that
$p_n: E_n \ra B_n$ is a weak homology fibration for any $n \geq 0$.
Assume that any face map $d_i: [n] \ra [n+1]$
induces an integral homology equivalence on homotopy fibres $Fib_{v}(p_{n+1}) \ra Fib_{d_i v}(p_n)$
for any vertex $v \in B_{n+1}$.
Then $p: || E_\bullet || \ra || B_\bullet ||$ is a homology fibration as well.
Moreover, if $w$ is a vertex in $|| B_\bullet ||$
lying in the same connected component as a vertex $v \in B_n$,
then there is a homology equivalence $Fib_w(p) \ra Fib_{v} (p_n)$.
\end{theorem}

\begin{proof}
Each step is a homotopy push-out involving only the face maps, so
Proposition~\ref{pushouthomologyfibration} applies. Hence $|| p_\bullet ||_n$ is a homology
fibration for any $n \geq 0$ and we conclude by Proposition~\ref{telescopehomologyfibration}.
\end{proof}

\medskip

One could actually prove a more general statement involving a colimit
over a small indexing category instead of the realization of a simplicial
space. In this paper we will not need such a statement.

\section{The generalized group completion}
\label{completion}

The aim is to find a model for the loops on the classifying space of a simplicial
category.
Let us start with a brief reminder on simplicial categories. More details can
be found for example in \cite[Section~1]{MR99k:57036}, especially about the link
with $2$-categories. Roughly speaking a simplicial category is a category equipped
with spaces of morphisms instead of sets of morphisms. So $mor_\M(i, j)$ is a space
for any objects $i, j \in \M$ and $mor_\M(i, i)$ contains the identity morphism as
distinguished base point. More precisely a simplicial
category $\M$ is a simplicial object in the category of small categories with constant
object set. It is helpful to look at $\M$ as a functor $\Delta^{op} \ra CAT$, where
the category of $n$-simplices is the category having same objects as $\M$ and morphisms
from $i$ to $j$ are the $n$-simplices of the space of morphisms from $i$ to $j$.
Taking now the nerve of this simplicial category degree by degree produces
a simplicial space denoted by $B\M_\bullet$, the {\em simplicial classifying space}.

A contravariant diagram $F: \M^{op} \ra Spaces$ is
the data of spaces $F(i)$ for all objects $i \in \M$ and natural continuous maps
$\mu_{i, j}: mor_\M(i, j) \times F(j) \ra F(i)$. The simplicial category itself
produces an example of diagram with $\M(i) = \coprod_{j \in Obj(\M)} mor_\M(i, j)$.

\begin{definition}
{\rm The {\em bisimplicial Borel construction} of a diagram $F: \M^{op} \ra Spaces$ is
the simplicial space $E_\M F_\bullet$ whose space of $n$-simplices is the disjoint
union over all $n$-tuples of objects in $\M$
$$
\coprod_{i_0, \dots, i_n} mor_\M(i_n, i_{n-1}) \times \cdots \times mor_\M(i_1, i_{0}) \times F(i_0)
$$
The degeneracy maps are the obvious inclusions. The face map
$d_n: E_\M F_n \ra E_\M F_{n-1}$ is projection on the last $n$ factors,
$d_0 = 1 \times \mu_{i_1, i_0}$, and the other $d_k$'s are defined by
composition $mor_\M(i_{k+1}, i_{k}) \times  mor_\M(i_k, i_{k-1}) \ra mor_\M(i_{k+1}, i_{k-1})$
.}
\end{definition}

The trivial diagram $T(i) = \{i\}$ is the diagram in which any morphism
$i \ra j$ induces the unique map $\{j\}\ra \{i\}$.
The bisimplicial Borel construction of the trivial diagram is nothing but
the simplicial classifying space of $\M$, i.e. $E_\M T_\bullet = B \M_\bullet$. Every
diagram $F: \M^{op} \ra Spaces$ comes with a natural transformation $\pi: F \ra T$ and hence we
get a map of simplicial spaces
$$
E_\M \pi_\bullet: E_\M F_\bullet \ra B\M_\bullet.
$$
The preimage of $\{i\}$ in the bisimplicial Borel construction is $F(i)$. Denote
by $E_\M F$ the realization $|| E_\M F_\bullet ||$, by $B\M$ the realization $|| B\M_\bullet ||$,
and by $\pi_\M: E_\M F \ra B\M$ the map induced by $\pi$.
We are ready to prove now the main theorem.

\begin{theorem}{\rm \cite[Theorem 3.2]{MR99k:57036}}
\label{groupcompletion}
Let $\M$ be a simplicial category and $F: \M^{op} \ra Spaces$ a
contravariant diagram. Assume that any morphism $f: i \ra j$
induces an isomorphism in integral homology $H_*(F(j); \Z) \ra
H_*(F(i); \Z)$. Then, for each object $i \in \M$, the map $F(i)
\ra Fib_i(\pi_\M)$ to the homotopy fibre of $\pi_\M$ over $i$ is a
homology equivalence.
\end{theorem}

\begin{proof}
We apply Theorem~\ref{simplicialhomologyfibration} to the map $E_\M \pi_\bullet$.
For any $n \geq 0$, the map $E_\M F_n \ra B\M_n$ is the projection on the first
factors, thus a (homology) fibration. As all faces but $d_0$ induce
the identity on the fibres, we have only to check that the face map $d_0$
induces a homology equivalence on the fibres. Choose a vertex
$$
(f_n, \, \dots \, , f_1, i_0)
\in mor_\M(i_n, i_{n-1}) \times \cdots \times mor_\M(i_1, i_{0}) \times \{i_0\}
$$
Its zeroth face is $(f_n, \, \dots \, , f_2, i_1)$ and the map induced on the
homotopy fibres is $F(f_0): F(i_0) \ra F(i_1)$. This is a homology equivalence
by assumption and we are done.
\end{proof}

\medskip

In order to identify the space $\Omega B \M$ we need to find
a diagram $F$ which satisfies the assumptions of Theorem~\ref{groupcompletion} and
for which the bisimplicial Borel construction $E_\M F$ is contractible.
We give a partial answer to that question which covers the applications made
in the context of the mapping class group.

Let us consider for any object $j \in \M$ the diagram $\M_j$ as defined
in~\cite[Section 3]{MR99k:57036}. It is the restriction of the diagram $\M$, i.e.
$\M_j(i) = mor_\M(i, j)$. This diagram
has a contractible bisimplicial Borel construction $E_\M \M_j \simeq *$
(see \cite[Lemma 3.3]{MR99k:57036}. Now fix an object $1 \in \M$ and
an endomorphism $\alpha: 1 \ra 1$, i.e. a vertex in the space of morphisms
$mor_\M(1, 1)$. Form the diagram $\M_\infty (i) =
hocolim ( \M_1(i) \stackrel{\alpha_*}{\lra} \M_1(i) \stackrel{\alpha_*}{\lra} \dots )$.
Since homotopy colimits commute with themselves $E_\M {\M_\infty} \simeq hocolim E_\M \M_1$ is
contractible and the homotopy fibre of $\pi_\M$ is $\Omega B\M$. We apply now the theorem
to the diagram $\M_\infty$.

\begin{proposition}
\label{burp}
Let $\M$ be a simplicial category and assume that there exists an endomorphism $\alpha$ of
a specific object $1$ such that any morphism $f: i \ra j$ induces an integral homology
equivalence $\M_\infty (j) \ra \M_\infty (i)$. Then the natural map
$\M_\infty (i) \ra \Omega B\M$ is an integral homology equivalence for any object
$i \in \M$. \hfill{\qed}
\end{proposition}

Finally one particularly likes the case when $\Omega B\M$ can be identified as Quillen's
plus construction applied to the space $\M_\infty (1)$. This means that the map
$\M_\infty (1) \ra \Omega B\M$ is not only a homology equivalence, but an acyclic map
(its homotopy fibre is acyclic). When is this so? In general a homology equivalence
is acyclic if the fundamental group of the base space acts nilpotently on the homology
of the homotopy fibre (assuming the fibre is connected, see \cite[4.3 (xii)]{MR84g:18028}).
This is usually rather difficult to verify. A stronger condition is that
$\pi_1 \M_\infty (1)$ is perfect. Then indeed every component of $\M_\infty (1)^+$
is 1-connected and hence $\M_\infty (1)^+$ is an $H\Z$-local space. Consider now the following
commutative square in which all arrows are homology equivalences
\[\xymatrix{
\M_\infty (1) \rto \dto & \Omega B\M \dto \\
\M_\infty (1)^+ \rto & (\Omega B\M)^+ \\
}\]
First $(\Omega B\M)^+ \simeq \Omega B\M$ since the fundamental group of any component
of a loop space is abelian. Moreover a loop space is always $H\Z$-local, so that
$\M_\infty (1)^+ \ra \Omega B\M$ is a homology equivalence between $H\Z$-local
spaces, thus a homotopy equivalence.

The above condition on the diagram $\M_\infty$ are precisely those checked in the proof of
\cite[Theorem~3.1]{MR99k:57036} to identify the plus construction on the classifying space
of the stable mapping class group as a loop space, which turns then out to be an infinite loop
space.

\begin{remark}
\label{conclude}
{\rm The homology theory which has been used in the present work
is integral homology and all applications we know of are obtained
working with integral homology. However, with little effort one
can replace this homology theory by an arbitrary (possibly
extraordinary) homology theory $E_*$. Hence an {\it
$E_*$-fibration} is a map $p: E \ra B$ such that $dp(\sigma) \ra
Fib_\sigma(p)$ is an $E_*$-equivalence. This is equivalent to
require that $p$ be a {\it weak $E_*$-fibration}, i.e. $dp(\sigma)
\ra dp(\theta \sigma)$ is an $E_*$-equivalence for any simplex
$\sigma$ in $B$ and any simplicial operation $\theta$. Then one
can prove the analogous of
Theorem~\ref{simplicialhomologyfibration}: The realization of a
natural transformation $p_\bullet: E_\bullet \ra B_\bullet$ of
simplicial spaces where all fibers have the same $E_*$-homology is
an $E_*$-fibration. The generalized group completion theorem has
an $E_*$-analogue as well, and the question would then be to
compare the homotopy type of $\Omega B \M$ with the
$E_*$-theoretic plus construction.}
\end{remark}

\section{Simplices versus topology}
\label{simplicesvstopology}

The general idea behind simplicial sets is to replace topological
data (points) by a combinatorial one (simplices). This is
precisely why one defines simplicially a homology fibration by
imposing a condition on the preimages of simplices, instead of
classically looking at preimages of points. There is however a
subtle difference, as shown by the following example due to
W.~Waldhausen, which we learned from J.~Rognes during the BCAT02.
A {\it simple map} of topological spaces is a map $f: X \ra Y$
such that the preimages of points $f^{-1}(y) \simeq *$ are
contractible for all $y \in Y$. Thus one would be tempted to
define simplicially a simple map as a map of spaces $f:X \ra Y$
for which preimages of simplices $dp(\sigma)\simeq *$ are all
contractible. This is not equivalent to the topological
definition. Consider indeed your favorite (but non-trivial)
acyclic space $A$. The map $A \ra *$ induces one on the unreduced
suspensions $\Sigma A \ra \Delta[1]$. The preimage of the
simplices in $\Delta[1]$ are either points, or $\Sigma A$, so all
are contractible. But topologically the geometric realization of
this map is not simple because the preimage of any other point
than the end points of the interval is $A$.

Recall that a map of topological spaces is a homology fibration if the
preimages of all points have the same homology type as the homotopy
fibre of $p$. We prove in this section that the simplicial and
topological definitions of homology fibrations are equivalent.
Basically this is due to the Mayer--Vietoris Theorem. The idea is
to take the barycentric subdivision of the map and reconstruct the
preimage of the barycenter of a simplex in the base from the data
given by the preimages of the simplices. Let us first recall some
standard definitions from \cite{MR19:759e} (or \cite[Chapter~4]{MR92d:55001}).

\medskip

Let $\mu$ be a proper face of $\Delta[n]$. We denote by $k_\mu$
the dimension of $\mu$, that is $\mu$ is  an injection $\mu~:
\Delta[k_\mu] \hookrightarrow \Delta[n]$. The {\it barycentric subdivision} of
$\Delta[n]$, denoted by $\Delta'[n]$, is the space which has as
$q$-simplices $\mathbf{\mu}$ the increasing sequences of $q+1$
faces of $\Delta[n]$, i.e. $\mu = (\mu_0, \cdots, \mu_q)$
where $ \mu_i(\Delta[k_i]) \subset
\mu_{i+1}(\Delta[k_{i+1}])$ for all $i \leq q-1$. The simplicial operations are the
usual: If $\theta~:\Delta[q] \lra \Delta[p]$  is any simplicial operation
then $\Delta' \alpha (\mu) =(\mu_{\theta(0)}, \cdots,
\mu_{\theta(q)})$.

The {\it subdivision} functor $Sd$ is left adjoint to Kan's
extension functor $Ex$ (see \cite[Section~7]{MR19:759e}). For any
space $E$, the $q$-simplices of $SdE$ are by definition the
equivalence classes $[x,\mu]$ of a simplex $x \in E$ of dimension
$p$ and $\mu \in \Delta'[p]$ of dimension $q$. Two pairs $(x,
\mu)$ and $(x', \mu')$ are equivalent if there exists a map
$\alpha: \Delta[p'] \ra \Delta[p]$ such that $x' = x \alpha$ and
$\mu = \Delta' \alpha (\mu)$. In other words, $SdE$ is the colimit
over the simplex category of $E$ of the subdivisions of these
simplices: $SdE = colim_{\Delta E} \Delta[n]'$.

Let us fix a surjective map $f:E \ra \Delta[n]$. Its subdivision
$Sdf: SdE \ra \Delta'[n]$ is defined as follows. Let $[x,\mu]$ be
a simplex in $SdE$ as above and consider for any $0 \leq i \leq q$
the composite
$$
\Delta[k_i] \stackrel{\mu_i}{\lra} \Delta[p] \stackrel{x}{\lra} E \stackrel{f}{\lra} \Delta[n]
$$
It can be decomposed in a unique way as a degeneracy followed by
an injection $\Delta[k_i] \stackrel{\phi_i}{\lra} \Delta[l_i]
\stackrel{\nu_i}{\mono} \Delta[n]$. Set $f([x,\mu]) = \nu =
(\nu_0, \, \dots , \, \nu_q)$.

\begin{definition}
{\rm In $\Delta'[n]$ fix a vertex $\alpha$, i.e. a proper face of
$\Delta[n]$. The {\it star} of $\alpha$, $St(\alpha)$
is the subspace of $\Delta'[n]$ which has as simplices the sequences
$(\mu_0, \cdots, \mu_p)$ such that $\forall i \leq p$,
$\im \mu_i \supset \im \alpha$.
We will further denote by $ESt(\alpha)$ the preimage of $St(\alpha)$ under
$Sdf$.}
\end{definition}

\begin{lemma}
\label{star}
The inclusion $Sdf^{-1}(\alpha) \mono ESt(\alpha)$ is a homotopy equivalence.
\end{lemma}

\begin{proof} Let $\alpha$ be of dimension $k$.
We construct first a retraction $r: ESt(\alpha) \ra Sdf^{-1}(\alpha)$.
Let $[x, \mu] \in ESt(\alpha)$ be a simplex of dimension~$q$. Then, for any $ i \leq q$,
there exists a maximal injective morphism $\Delta[t_i] \mono \Delta[k_i]$ (determined by
the vertices of $\mu_i$ whose image under $f(x)$ is a vertex of $\alpha$) together
with a (necessary unique) surjection $\phi~: \Delta[t_i] \ra \Delta[k]$
rendering the following diagram commutative
\[
\xymatrix{
\Delta[k_i] \rmono^{\mu_i} & \Delta[p]  \rto^{x} &  E \dto^f \\
\Delta[t_i] \ar@{^{(}->}[u] \rto^{\phi} & \Delta[k] \rto^{\alpha} & {\Delta[n]}
}
\]
We denote the composite $\Delta[t_i] \ra \Delta[k_i] \ra \Delta[p]$ by $\bar\mu_i$
and define $r[x, \mu] = [x, \bar\mu]$. By construction $Sdf([x, \bar\mu])$ is some
degeneracy of $\alpha$. Moreover $r$ is well defined and is clearly a retraction of
the inclusion $i: Sdf^{-1}(\alpha) \mono ESt(\alpha)$.

Finally we construct a homotopy $H: ESt(\alpha) \times \Delta[1]
\ra ESt(\alpha)$ from $i \circ r$ to the identity. Let $([x, \mu],
\tau)$ be a $q$-simplex in the cylinder, so $\tau$ is a
$q$-simplex in $\Delta[1]$ and can be represented by a sequence of
$r+1$ zero's and $q-r$ one's: $(0 \dots 0 1 \cdots 1)$. Define
then $H([x, \mu], \tau) = [x, \bar\mu_0, \dots , \bar\mu_r,
\mu_{r+1}, \dots, \mu_q]$.
\end{proof}

\bigskip

In the next proposition we use the decomposition of $\Delta'[n]$
as union of all its stars. More precisely consider the category
${\mathcal{C}}_n$ whose objects are the non-degenerate simplices of
$\Delta[n]$ and whose morphisms are generated by the faces $\sigma
\ra d_i\sigma$. The unique non-degenerate simplex $\tau$ of
dimension $n$ is an initial object and diagrams indexed by
 ${\mathcal{C}}_n$ are $n$-cubes without terminal object. We have $\Delta'[n] =
colim_{\sigma \in {\mathcal{C}}_n} St(\sigma) = hocolim_{\sigma \in
{\mathcal{C}}_n} St(\sigma)$ because the diagram $St$ is cofibrant (see
for example \cite{wgdspal}), and even strongly co-Cartesian as
defined in \cite[Definition~2.1]{MR93i:55015}. Likewise
$$
E \simeq SdE = colim_{\sigma \in {\mathcal{C}}_n} ESt(\sigma) =
hocolim_{\sigma \in {\mathcal{C}}_n} ESt(\sigma)
$$

\begin{proposition}
\label{barycenter}
Let $f: E \ra \Delta[n]$ be a homology fibration. Then the
preimage of the barycenter of $\Delta'[n]$ under $Sdf$ has the
same homology type as $E$. In particular the realization $|f|: |E|
\ra |\Delta[n]|$ is a homology fibration of topological spaces.
\end{proposition}

\begin{proof}
By Lemma~\ref{star} the values of the cubical diagram $ESt$ are
equivalent to the preimages $Sdf^{-1}(\sigma)$. When $\sigma$ is a
vertex of $\Delta[n]$, one has that $Sdf^{-1}(\sigma) \simeq
f^{-1}(\sigma) = df(\sigma)$, which by hypothesis has the same
homology type as $E$. By induction on the dimension of $\sigma$ we
can assume thus that all values in the diagram but the initial one
($ESt(\tau) \simeq Sdf^{-1}(\tau)$, the preimage of the
barycenter) are homology equivalent to $E$. As the homotopy
colimit of the cubical diagram is $E$, we deduce that $ESt(\tau)$
as well has the same homology type as $E$. We claim that this
implies that $|f|$ is a (topological) homology fibration. Indeed
by induction again we need only to compute preimages under $|f|$
of points in the interior of the realization of $\Delta[n]$. Any
such preimage is a deformation retract of the preimage under $|p|$
of the open simplex, so it is enough to consider the barycenter.
The above computation shows precisely that it has the same
homology type as $|E|$, the homotopy fibre of~$|f|$.
\end{proof}

\bigskip

Let us now consider a map $p: E \ra B$. To compare both types of
homology fibrations we need to control the homological properties
of fibers of points in the realization of spaces. Any point $b \in
|B|$ lies in the interior of the realization of a unique
non-degenerate simplex $\sigma_b \in B$ (see for instance
\cite[Lemma 4.2.5]{MR92d:55001}). Moreover the interior of the
realization of $\sigma_b$ embeds in $|B|$.

\begin{theorem}
A map of spaces $p: E \ra B$ is a homology fibration if and only
if its realization $|p|: |E| \ra |B|$ is a homology fibration of
topological spaces.
\end{theorem}

\begin{proof}
First assume that $p: E \rightarrow B$ is a homology fibration. We
need to compute the homology type of fibers of points in the
realization of $B$ and show that the map $|p|^{-1}(b) \rightarrow
Fib_b (|p|)$ is a homology equivalence , where $Fib_b (|p|)$
denotes the homotopy fiber of $|p|$ over the connected component
of $b$. When $\sigma = \sigma_b$ is a $0$-simplex, this is trivial
as $p$ is a homology fibration. If $\sigma$ is of dimension $n
\geq 1$, notice that all the fibers over the points in the
interior of $|\sigma|$ have the same homotopy type (a
straightforward computation shows then that the preimage any point
is a deformation retract of the preimage under $|p|$ of the open
simplex). Therefore it suffices to analyze the barycenter
$\iota_{n}$ of the realization of $\sigma$ and to prove that
$|p|^{-1}(\iota_{n}) \rightarrow Fib_{\iota_{n}} (|p|)$ is a
homology equivalence. As the realization functor commutes with
finite limits (see \cite[Theorem 4.3.16]{MR92d:55001}), we have a
pull-back square~:
\[\xymatrix{
{| dp(\sigma) |} \dto  \rto  & {|E|}\dto\\
{| \Delta[n] |} \rto^{\sigma} & {|B|} }
\]
The map $dp(\sigma) \ra \Delta[n]$ is a homology fibration as the
pull-back of any simplex of the base $\Delta[n]$ coincides with
the pull-back of a simplex of $B$ along $p$, which has the same
homology type as $dp(\sigma)$. By Proposition~\ref{barycenter},
its realization is a homology fibration: The preimage of the
barycenter of $| \Delta[n] |$ is homology equivalent to the
homotopy fibre $|dp(\sigma)|$, which by assumption has the same
homology type as the homotopy fibre $|F|$ of $|p|$.

Assume now $|p|: |E| \rightarrow |B|$ is a homology equivalence.
Inductively we may suppose that for all simplices of dimension
$\leq n-1$ the pull-back $dp(\tau)$ is homology equivalent to the
homotopy fibre above the component of $\tau$. Let $\sigma$ be a
simplex of dimension $n$. We have as before a pull-back diagram
\[\xymatrix{
{| dp(\sigma) |} \dto  \rto  & {|E|}\dto\\
{| \Delta[n] |} \rto^\tau & {|B|} }
\]

Decompose $dp(\sigma)$ as a cubical homotopy colimit
$dp(\sigma) \simeq hocolim_{\tau \in \mathcal{C}_n} ESt(\tau)$
following the method seen in the proof of
Proposition~\ref{barycenter}. As $|p|$ is a homology fibration,
there is a natural transformation by homology equivalences to the
constant cubical diagram $Fib_\sigma (p)$ (use Lemma~\ref{star}).
A homotopy colimit of homology equivalences is a homology
equivalence, hence $dp(\sigma) \rightarrow Fib_\sigma (p)$ is a
homology equivalence as well.
\end{proof}

\bibliographystyle{alpha}\label{bibliography}

\begin{thebibliography}{99}

\bibitem{MR80d:55001} Adams, J. F. \it {Infinite loop spaces}, volume 90. Princeton University Press, 1978. Annals of Mathematics Studies

\bibitem{Amit} Amit, A.  Direct limits over categories with contractible nerve, \emph{Master Thesis, The Hebrew University of Jerusalem}, 1994.

\bibitem{MR84g:18028} Berrick, A. J. \it{An approach to algebraic ${K}$-theory}, Pitman (Advanced Publishing Program), Boston, Mass., 1982.
 
 \bibitem{MR47:3489} Baratt, M. and Priddy, S. On the homology of non-connected monoids and theirs associated groups, \it{Comment. Math. Helv.}, 47:1-14, 1972.
 
 \bibitem {ChSc} Chach\'olski, W. and Scherer, J. Homotopy Theory of Diagrams, \it {Mem. Amer. Math. Soc.},  155 (736):  ix+90, 2002.
 
 \bibitem{dror:book} Dror Farjoun, E. \it{Cellular spaces, null spaces and homotopy localization}, volume 1622 of  \it{Lecture Notes in Mathematics}, Springer-Verlag, Berlin, 1996.
 
 \bibitem{wgdspal} Dwyer, W. G. and Spali{\'n}ski, J. Homotopy theories and model categories. In Handbook of algebraic topology, pages 73--126. North-Holland, Amsterdam, 1995.

 \bibitem{MR95a:14023} Friedlander, E.M. and Mazur, B. Filtrations on the homology of algebraic varieties. \it{Mem. Amer. Math. Soc.}, 110 (529) : x + 110, 1994. With an appendix by Daniel Quillen.

\bibitem{MR92d:55001} Fritsch, R. and Piccinini, R. A. \it{Cellular structures in topology}, volume 19 of \it{Cambridge Studies in Advanced Mathematics}. Cambridge University Press, Cambridge, 1990.

\bibitem{MR93i:55015} Goodwillie, T. G.  Calculus. {I}{I}. {A}nalytic functors, \it {$K$-Theory}, 5 (4) : 295--332, 1991/92.

\bibitem{MR90j:55029} Jardine, J. F. The homotopical foundations of algebraic ${K}$-theory. In \it{Algebraic $K$-theory and algebraic number theory (Honolulu,
              HI, 1987)}, pages 57--82. Amer. Math. Soc., Providence, RI, 1989.
 
 \bibitem{MR19:759e} Kan, D. M. On c.s.s. complexes, \it{Amer. J. Math.}, 79 : 449--476, 1957.

\bibitem{MR51:6806} May, J. P.  Classifying spaces and fibrations, \it{Mem. Amer. Math. Soc.}, 1 (1, 155) :  xiii+98, 1975.

\bibitem{MR91j:18021} Moerdijk, I. Bisimplicial sets and the group-completion theorem. In \it{Algebraic $K$-theory: connections with geometry and topology (Lake Louise, AB, 1987)}, pages 
225--240, Kluwer Acad. Publ., Dordrecht, 1989.

\bibitem{MR53:6547} McDuff, D. and Segal, G. Homology fibrations and the ``group-completion'' theorem. \it{Invent. Math.}, 31 (3) : 279--284, 1975/76.

\bibitem{MR43:5525} Quillen, D. The Adams conjecture, \it{Topology}, 10 : 76--80, 1971.

\bibitem{MR50:5782} Segal, G. Categories and cohomology theories. \it{Topology}, 13 : 293--312, 1974.

\bibitem{MR99k:57036} Tillmann, U. On the homotopy of the stable mapping class group. \it{Invent. Math.}, 130 (2) : 257--275, 1997.
\end{thebibliography}

\medskip

\noindent
Wolfgang Pitsch, Institut de Math\'ematiques, Universit\'e de Gen\`eve,
CH--1200 Gen\`eve, Switzerland, e-mail: {\tt wpitsch@ima.unige.ch}

\medskip

\noindent
J\'er\^ome Scherer, Departament de Matem\'atiques, Universitat Aut\`onoma de
Barcelona, E--08193 Bellaterra, Spain, e-mail: {\tt
jscherer@mat.uab.es}

\end{document}